\documentclass[12pt,oneside,reqno]{amsart}
\usepackage{}
\usepackage{amssymb}
\usepackage{bbm}
\usepackage{cases}
\usepackage{amsmath}
\usepackage{graphicx}
\usepackage{cite}
\usepackage{mathrsfs}
\usepackage{stmaryrd}
\usepackage{amsfonts}
\usepackage{enumerate,amsmath,amssymb,amsthm}

\numberwithin{equation}{section}

\newcommand{\be}{\begin{eqnarray}}
\newcommand{\ee}{\end{eqnarray}}
\newcommand{\ce}{\begin{eqnarray*}}
\newcommand{\de}{\end{eqnarray*}}
\newtheorem{theorem}{Theorem}[section]
\newtheorem{lemma}[theorem]{Lemma}
\newtheorem{remark}[theorem]{Remark}
\newtheorem{definition}[theorem]{Definition}
\newtheorem{proposition}[theorem]{Proposition}
\newtheorem{Examples}[theorem]{Example}
\newtheorem{corollary}[theorem]{Corollary}

\def\e{{\mathrm{e}}}
\def\eps{\varepsilon}

\def\[{{\Big[}}
\def\]{{\Big]}}
\def\<{{\langle}}
\def\>{{\rangle}}
\def\({{\Big(}}
\def\){{\Big)}}

\def\bx{{\mathbf{x}}}

\def\dif{{\mathord{{\rm d}}}}

\def\no{\nonumber}
\def\={&\!\!=\!\!&}
\def\bt{\begin{theorem}}
\def\et{\end{theorem}}
\def\bl{\begin{lemma}}
\def\el{\end{lemma}}
\def\br{\begin{remark}}
\def\er{\end{remark}}

\def\bd{\begin{definition}}
\def\ed{\end{definition}}
\def\bp{\begin{proposition}}
\def\ep{\end{proposition}}
\def\bc{\begin{corollary}}
\def\ec{\end{corollary}}
\def\bx{\begin{Examples}}
\def\ex{\end{Examples}}

\def\cB{{\mathcal B}}
\def\cC{{\mathcal C}}

\def\cH{{\mathcal H}}

\def\cT{{\mathcal T}}

\def\mE{{\mathbb E}}

\def\mI{{\mathbb I}}

\def\mN{{\mathbb N}}

\def\mP{{\mathbb P}}
\def\mQ{{\mathbb Q}}
\def\mR{{\mathbb R}}

\def\sF{{\mathscr F}}

\def\sL{{\mathscr L}}

\def\geq{\geqslant}
\def\leq{\leqslant}

\allowdisplaybreaks

\begin{document}

\title{Well-posedness and long time behavior of singular Langevin stochastic differential equations}

\date{}

\author{Renming Song\ \ and \ \ Longjie Xie}

\address{Renming Song:
Department of Mathematics, University of Illinois,
Urbana, IL 61801, USA\\
Email: rsong@illinois.edu
 }

\address{Longjie Xie:
School of Mathematics and Statistics, Jiangsu Normal University,
Xuzhou, Jiangsu 221000, P.R.China\\
Email: xlj.98@whu.edu.cn
 }

\thanks{Research of R. Song is supported by the Simons Foundation (\#429343, Renming Song). L. Xie is supported by NNSF of China (No. 11701233) and NSF of Jiangsu (No. BK20170226). The Project Funded by the PAPD of Jiangsu Higher Education Institutions is also gratefully acknowledged}

\begin{abstract}
In this paper, we study damped Langevin stochastic differential equations with singular velocity fields. We prove the strong well-posedness of such equations.  Moreover, by combining the technique of Lyapunov functions with Krylov's estimate, we also establish the exponential ergodicity for the unique strong solution.

\bigskip

\noindent\textbf{AMS 2010 Mathematics Subject Classification:} Primary 60H10;
Secondary 37A25; 82C31

\medskip

  \noindent\textbf{Keywords and Phrases:} Pathwise uniqueness; Langevin equation; Krylov's estimate; exponential ergodicity.
\end{abstract}

\maketitle \rm

\section{Introduction and main results}

In this paper, we study the following Langevin stochastic differential equation (SDE):
\begin{equation} \label{sde1}
\left\{ \begin{aligned}
&\dif X_t=V_t\dif t,\qquad\qquad\qquad\qquad\qquad\qquad\qquad\quad\,\, X_0=x\in\mR^{d},\\
&\dif V_t=-\gamma V_t\dif t-\nabla F(X_t)\dif t+G(V_t)\dif t+\dif W_t,\quad V_0=v\in\mR^{d},
\end{aligned} \right.
\end{equation}
where $W_t$ is a standard $d$-dimensional  Brownian motion defined on some filtered probability space $(\Omega,\sF,(\sF_t)_{t\geq 0},\mP)$, and $\gamma>0$ is a friction constant which ensures that our Hamiltonian is damped.
The coefficient $F: \mR^d\rightarrow\mR$ is a non-negative smooth function,
and $G: \mR^{d}\rightarrow\mR^{d}$ is a Borel function whose assumption will be given later.
The equation \eqref{sde1} describes the position and momentum of a particle of unit mass that moves according to Newton's second law in a smooth potential $F$
that is subject to friction and noise. Such equation appears in various applications such as stochastic non-linear oscillators, surface diffusions and polymer dynamics, see \cite{Du,Ma-We,Pa}.

\vspace{1mm}
Note that \eqref{sde1} is degenerate with dissipation and noise appearing only in
the momentum part.
However, there will be sufficient interaction between the position and momentum
so that dissipation and noise get transmitted from the equation for $V_t$ to the equation for $X_t$, which ultimately leads to the ergodicity and exponential convergence to equilibrium. In the case $G\equiv0$, it is well known that \eqref{sde1} is well-posed and many people have studied the asymptotic behavior of the unique strong solution, see \cite{B-S,C-F,D-M-S,E,H-M,H-M2,O-P-P,Ma,M-S-H,Re,T,Wu} and the references therein.
Most of these papers rely on the existence of a Lyapunov function,
which means that the force should be strong enough to make the system return quickly to compact subsets of $\mR^{2d}$. In this case, the generator of the SDE \eqref{sde1} is a hypoelliptic operator given by
$$
\sL:=v\cdot\nabla_x-\gamma v\cdot\nabla_v-\nabla F(x)\cdot\nabla_v+\tfrac{1}{2}\Delta_v=:\sL_1+\sL_2,
$$
where
$$
\sL_1:=\tfrac{1}{2}\Delta_v-\gamma v\cdot\nabla_v
$$
is the generator of an Ornstein-Uhlenbeck process, and
$$
\sL_2:=v\cdot\nabla_x-\nabla F(x)\cdot\nabla_v
$$
is the Liouville operator generated by the Hamiltonian $H(x,v):=\tfrac{1}{2}|v|^2+F(x)$.
Under the condition that there exist constants $C_0>0$ and $\vartheta\in (0,1)$ such that
\begin{align}
\frac{1}{2}\<x,\nabla F(x)\>\geq \vartheta F(x)+\gamma^2\frac{\vartheta(2-\vartheta)}{8(1-\vartheta)}|x|^2-C_0,  \label{FF}
\end{align}
one can check that the following Hamiltonian functional with lower order perturbation
\begin{align}
\cH(x,v):=H(x,v)+\tfrac{\gamma}{2}\<x,v\>+\tfrac{\gamma^2}{4}|x|^2+1 \label{Lya}
\end{align}
is a good choice of a Lyapunov function (see \cite{M-S-H}), i.e., for some positive constants $c_0$, $K_1$ and all $x,v\in\mR^d$,
\begin{align}
\cH(x,v)\geq 1+\tfrac{\gamma^2}{12}|x|^2+\tfrac{1}{8}|v|^2\quad\text{and}\quad\sL\cH(x,v)\leq -c_0\cH(x,v)+K_1. \label{lay}
\end{align}
Recently, the case of the SDE \eqref{sde1} with a singular potential field $F$ has been studied in \cite{C-G,C-H-M-M-S,G-S}. The exponential ergodicity  was also obtained by constructing explicit Lyapunov functions, see also the recent work \cite{H-Ma} and the references therein.

\vspace{1mm}
We intend to study the existence and uniqueness of strong solutions as well as the exponential ergodicity  of \eqref{sde1} under the presence of a singular velocity field $G$, which destroys the dissipation in the momentum part and makes the classical Lyapunov
condition very difficult to check, if possible at all.
Let us first recall the following well-known concept, see for instance, \cite{Go-Ma} and \cite{Wu}.

In the remainder of this paper, $\cH$ always stands for the function defined in \eqref{Lya}.

\bd
The invariant distribution $\mu$ (if exists) of an $\mR^{2d}$-valued Markov process
$M_t$ is said to be $\cH$-uniformly exponentially ergodic, if there exist constants $C,\eta>0$ such that for
all $y\in\mR^{2d}$ and all Borel functions $f: \mR^{2d}\rightarrow\mR^d$ with
$|f|\leq \cH$,
\begin{align*}
\big|\mE^{y}f(M_t)-\mu(f)\big|\leq C\cH(y)\e^{-\eta t},\quad \forall t\geq 0,
\end{align*}
where $\mE^y$ is the expectation with respect to
$\mP^y$, the law of $M$ with initial value  $M_0=y$,
and $\mu(f)$ denotes the integral of $f$ respect to $\mu$.
\ed

The main result of this paper is as follows.

\bt\label{main1}
Let $G\in L^p(\mR^d)$ with $(2\vee d)<p\leq \infty$.
For each $y:=(x,v)^{\mathrm{T}}\in\mR^{2d}$, \eqref{sde1} admits a unique strong solution $Y_t=(X_t,V_t)^\mathrm{T}$.
Moreover, if we further assume that \eqref{FF} holds  together with one of the following conditions:
\begin{enumerate}[{\bf(A)}]
\item $|\nabla F(x)|^2\leq C_1(1+|x|^2+F(x))$;

\item $|\nabla F(x)|\leq C_1(1+|x|^2+F(x))$ and $G\in L^p(\mR^d)$ with $2d<p\leq \infty$,
\end{enumerate}
where $C_1>0$ is a constant,  then $Y_t$ has a unique invariant
distribution $\mu$ which is $\cH$-uniformly exponentially ergodic.
\et

Note that condition ${\bf(A)}$ includes the case of a harmonic potential $F(x)=\tfrac{1}{2}|x|^2$,
while any polynomial $F$ which grows  at infinity like $|x|^{2\ell}$ for some positive integer $\ell$ satisfies condition ${\bf(B)}$.
The presence of the flexible singular term $G$
can be used to describe stochastic non-linear oscillators (see \cite{Ma-We}) as well as degenerate particle systems with singular velocity fields arising in mathematical physics (see \cite{K-T-Y,Y-T-K,S}). In particular, we have the following example.

\bx
Consider the following equation:
\begin{equation*}
\left\{ \begin{aligned}
&\dif X_t=V_t\dif t,\qquad\qquad\qquad\qquad\qquad\qquad\qquad\quad\,\, X_0=x\in\mR^{d},\\
&\dif V_t=-\gamma V_t\dif t-X_t\dif t+\tfrac{1}{|V_t|^\alpha}1_{\{|V_t|\leq K\}}\dif t+\dif W_t,\quad V_0=v\in\mR^{d},
\end{aligned} \right.
\end{equation*}
where $K>0$ and $0<\alpha<1$. Then all the conclusions in Theorem \ref{main1} hold.
\ex

It turns out that the main difficulty for
studying  \eqref{sde1} is that we have to treat simultaneously
the singular term $G$ and the super-linear growth part $F$ in the coefficients.
Here we give a brief description of the main idea of the proof and
make some comparisons with existing results.
To prove the strong well-posedness of \eqref{sde1}, we will use  Zvonkin's transform, see \cite{Zv}. We mention that stochastic kinetic equations with singular drifts have also been studied in \cite{F-F-P-V,Wa-Zh1,Zh}. However, their assumptions on the integrability of the coefficients are much stronger than ours and can not be applied to our case directly.
Compared with \cite{F-F-P-V,Kr-Ro,Wa-Zh1,Zh}, where global integrability of the coefficients are needed, the super-linear growth part $F$ 
will cause some trouble, particularly for deriving Krylov's estimate.
Since the uniqueness is a local concept, we will use a localization technique to truncate the growth coefficients.
Thus, no further conditions on $F$ are needed to ensure the uniqueness of \eqref{sde1}. 
In view of \cite{Kr-Ro}, our condition on $G$ should be the best possible.

To study the long time behavior of the unique strong solution,
the localization technique is of no help,
and this is why we need some growth conditions on $F$ to derive a global Krylov's estimate for the solution.
In particular, we need to solve the following quasi-linear elliptic equation in Sobolev spaces:
\begin{align*}
\lambda u-\tfrac{1}{2}\Delta u-G\cdot\nabla u-\kappa|\nabla u|^2=f,
\end{align*}
where $\lambda,\kappa\geq 0$, which is of independent interest. It is also interesting to note that conditions ${\bf(A)}$ and ${\bf(B)}$ reflect that some balances are needed between the integrability of $G$ and the growth property of $F$.
In the recent work \cite{XZ2}, the ergodicity of non-degenerate SDEs with singular dissipative coefficients has been studied. The idea in \cite{XZ2} is to use Zvonkin's transform to remove
the singular part of the drift and the fact that the dissipative property is preserved
by Zvonkin's transform in the setting of \cite{XZ2}.
However, such an idea can not be applied to our case. As we shall see, the Hamiltonian
structure of \eqref{sde1} will be totally destroyed by Zvonkin's transform.
Thus, it will be really difficult to find a Lyapunov function for the new equation (and thus for the original one).
To overcome this difficulty, we will use Krylov's
estimate to get a good control on the expectation of the singular part,
and then combine with the Lyapunov technique to get the
existence of invariant distributions for \eqref{sde1}. The uniqueness of invariant distribution follows from the strong Feller property and irreducibility of the unique strong solution.
Note that the argument used in \cite{Kr-Ro,XZ,XZ2,Zh1} to derive
the strong Feller property and irreducibility of the strong solution
does not apply any more since \eqref{sde1} is degenerate. Instead, we adopt the method developed in \cite{M-S} by making use of  the Girsanov transform.

\vspace{2mm}

The paper is organized as follows. In Section 2, we prove Krylov's estimate for the  solution to \eqref{sde1}. The proof of Theorem \ref{main1} will be given in Section 3.
Throughout this paper, we use the following convention: $C$ with or without subscripts will denote a positive constant, whose value may change from one appearance to another., and whose dependence on parameters can be traced from calculations.

\section{Krylov's estimate for the solution of \eqref{sde1}}

This section consists of three subsections. We first solve a quasi-linear elliptic equation in Sobolev spaces in Subsection 2.1. Then, we derive Krylov's estimates for the solution of \eqref{sde1} by using two different methods: the Girsanov transform method in Subsection 2.2 and the elliptic equation approach in Subsection 2.3.

\subsection{Quasi-linear elliptic equation}
Let us introduce some notation. For $1\leq p\leq \infty$, the norm in $L^p(\mR^d)$ will be denoted by $\|\cdot\|_p$, and for $k\in\mN$, $W^{k,p}(\mR^d)$ is the usual Sobolev space with norm
$$
\|f\|_{k,p}:=\|f\|_p+\Sigma_{|\alpha|=1}^k\|D^\alpha f\|_{p},
$$
where for any multi-index $\alpha=(\alpha_1, \dots, \alpha_d)$ with $|\alpha|=\alpha_1+\cdots+\alpha_d$,
$D^\alpha f$ stands for the weak derivative $\frac{\partial^{|\alpha|}f}{\partial x_1^{\alpha_1}\cdots\partial x_d^{\alpha_d}}$.
It is known that when $p>1$, a function $f\in L^p(\mR^d)$ belongs to $W^{1,p}(\mR^d)$ if and only if for a.e. $x,y\in\mR^d$, there  
exists $g\in L^p(\mR^d)$ such that
\begin{align}
|f(x)-f(y)|\leq |x-y|\big(g(x)+g(y)\big).   \label{w1p}
\end{align}
When $r>0$ is not an integer, the fractional Sobolev space
$W^{r,p}(\mR^d)$ is defined to be the space of functions with
\begin{align*}
\|f\|_{r,p}:=
\|f\|_{[r],p}+ \sum_{|\alpha|=k}\Bigg(\int_{\mR^d}\!\!\int_{\mR^d}\frac{|D^\alpha f(x)-D^\alpha f(y)|^p}{|x-y|^{d+(r-[r])p}}\dif x\dif y\Bigg)^{1/p}<\infty.
\end{align*}
The celebrated Sobolev's embedding theorem tells us that
if $0<r<d/p$, then for any $p\leq q\leq pd/(d-rp)$, it holds that
\begin{align}
W^{r,p}(\mR^d)\hookrightarrow L^q(\mR^d),\label{emb2}
\end{align}
and if $r>d/p$, we have
\begin{align}
W^{r,p}(\mR^d)\hookrightarrow C^{r-\frac{d}{p}}_b(\mR^d)\subseteq L^\infty(\mR^d),\label{emb}
\end{align}
where for $\gamma\geq 0$, $C_b^\gamma(\mR^d)$ denotes the usual H\"older space. Such kinds of embeddings will play important roles below.

To derive Krylov's estimate for the solution to \eqref{sde1}, we shall need to consider the following quasi-linear elliptic equation:
\begin{align}
\lambda u(x)-\tfrac{1}{2}\Delta u(x)-G(x)\cdot\nabla u(x)-\kappa|\nabla u(x)|^2=f(x),  \label{eq2}
\end{align}
where $\lambda,\kappa\geq 0$ are constants. To handle the non-linear part, we first show a basic result about a sequence of real numbers, which will be needed below.

Let $\delta>0$ be a constant. Define a sequence of 
positive numbers $\{\cC_n\}_{n\geq 1}$ recursively by
\begin{align}
\cC_1=1\quad\text{and}\quad \cC_{n+1}=1+\delta{\cC_{n}}^2.    \label{cn}
\end{align}
Note that $\cC_n$ is increasing. We have:

\bl\label{nu}
For every $0<\delta\leq 1/4$, the sequence $\{\cC_n\}_{n\geq 1}$
converges.
\el
\begin{proof}
Let $C_n$ be the $n^{\text{th}}$ Catalan number, that is, $C_0=1$ and for $n\geq 0$,
$$
C_{n+1}:=\sum_{i=0}^{n}C_{i}C_{n-i}.
$$
Note that
$$
\cC_1=C_1,\quad\text{and}\quad \cC_2=C_0+\delta C_1.
$$
By induction and the recurrence relation of $C_n$, we can show that
$$
\cC_n=\sum_{k=0}^{n-1}C_k\delta^k+o(\delta^n).
$$
The generating function of the Catalan numbers is
$$
\sum_{k=0}^{\infty}C_k\delta^k=\frac{1-\sqrt{1-4\delta}}{2\delta},
$$
which means that $\{\cC_n\}$ converges if $\delta\leq1/4$.
\end{proof}

Now let $p(t, x):=(2\pi t)^{-d/2}\exp(-\frac{|x|^2}{2t})$, which is the fundamental solution of $\frac{1}{2}\Delta$,
and denote by $\cT_t$ the corresponding semigroup, i.e.,
$$
\cT_tf(x):=\int_{\mR^d}f(y)p(t,x-y)\dif y,\quad \mbox{ for all }  f\in \cB_b(\mR^d).
$$
We now prove the following result.
The key point of the following proof is to use
Sobolev's embeddings to handle the non-linear term $|\nabla u|^2$.

\bl\label{ell}
Let $\kappa\geq 0$ be a constant.

\noindent{\bf(i)}
Suppose that  $G\in L^p(\mR^d)$ with $p\in (d, \infty)$. 
Then for every $f\in L^q(\mR^d)$ with $q\in (d/2, \infty)$, 
there exists a $\lambda_0>0$ such that for all $\lambda\geq\lambda_0$, \eqref{eq2} admits a unique solution $u\in W^{2,q}(\mR^d)$.
Moreover, we have
\begin{align}
\lambda\|u\|_q+\sqrt{\lambda}\|\nabla u\|_q+\|\nabla^2u\|_q\leq C_1\|f\|_q,  \label{ess}
\end{align}
where $C_1=C_1(d,\|G\|_p)>0$ is a constant.

\vspace{1mm}
\noindent{\bf(ii)}
Let $p$ and $q$ be as in part (i).
Given two sequences of functions $G_n,f_n\in C_0^\infty(\mR^d)$ such that $G_n\rightarrow G$ in $L^p(\mR^d)$ and $f_n\rightarrow f$ in $L^q(\mR^d)$. Let $u_n$ be the corresponding solution to \eqref{eq2} with $G, f$ replaced by $G_n, f_n$. Then we have $u_n\in C_b^2(\mR^d)\cap W^{2,q}(\mR^d)$ with
\begin{align}
\sup_{n\geq 1}\|u_n\|_{2,q}\leq C_2\quad\text{and}\quad\|u_n-u\|_{2,q}\leq C_2\big(\|f_n-f\|_q+\|G_n-G\|_p\big),    \label{ap}
\end{align}
where $C_2=C_2(d,\|G\|_p,\|f\|_q)$ is a positive constant.
\el

\begin{proof}
{\bf(i)} It is now standard that we only need to show the existence and uniqueness of solution to the integral equation
\begin{align}
u(x)=\int_0^\infty\!\e^{-\lambda t}\cT_t\big(f+G\cdot\nabla u+\kappa|\nabla u|^2\big)(x)\dif t,   \label{inte}
\end{align}
and prove the a priori estimate \eqref{ess}. Let us construct the solution via Picard's iteration.
For $\lambda>0$, set
\begin{align}
u_1(x):=\int_0^\infty\!\e^{-\lambda t}\cT_tf(x)\dif t. \label{u1}
\end{align}
By the classical theory of linear elliptic equations (see \cite[Lemma 4.2]{XZ2}, for instance), there exists a constant $C_d>0$ such that
\begin{align*}
\lambda\|u_1\|_q+\sqrt{\lambda}\|\nabla u_1\|_q+\|\nabla^2u_1\|_{q}\leq C_d\|f\|_{q}.
\end{align*}
For $n\geq 2$, define $u_n$ recursively by
\begin{align}
u_n(x):=\int_0^\infty\!\e^{-\lambda t}\cT_t\big(f+G\cdot\nabla u_{n-1}+\kappa|\nabla u_{n-1}|^2\big)(x)\dif t.  \label{ind}
\end{align}
We claim that there exists a $\delta>0$ such that for all $n\in\mN$,
\begin{align}
\lambda\|u_n\|_q+\sqrt{\lambda}\|\nabla u_n\|_q+\|\nabla^2u_n\|_{q}\leq C_d\,\cC_n\|f\|_{q},    \label{indu}
\end{align}
where the constant $\cC_n$ satisfies the recursion relation \eqref{cn}.

We first prove the above claim. It is clear that the assertion holds for $n=1$.
Assume that \eqref{indu} is true for some $n\geq 1$, then we can deduce
\begin{align*}
\lambda\|u_{n+1}\|_q+\sqrt{\lambda}\|\nabla u_{n+1}\|_q+\|\nabla^2u_{n+1}\|_{q}&\leq C_d\big(\|f\|_q+\|G\cdot\nabla u_n\|_q+\kappa\||\nabla u_n|^2\|_q\big).
\end{align*}
On one hand, recall that $G\in L^p(\mR^d)$ with $p\in (d, \infty)$.
By the Sobolev embedding \eqref{emb2}, it holds that
$$
\nabla u_n\in W^{1,q}(\mR^d)\subseteq W^{\alpha,q}(\mR^d)\hookrightarrow L^r(\mR^d) \quad\text{with}\quad \alpha=\tfrac{d}{p}<1\quad\text{and}\quad r=\tfrac{pq}{p-q}.
$$
Note that $1/q=1/p+1/r$, we can derive by Young's inequality that for every $\eps>0$, there exists a $C_\eps$ such that
\begin{align}
\|G\cdot\nabla u_n\|_q&\leq \|G\|_p\|\nabla u_n\|_r\leq \|G\|_p\|\nabla u\|_{\alpha,q}\no\\
&\leq \eps\|G\|_p\|\nabla^2u_n\|_{q}+C_\eps\|G\|_p\|u_n\|_q.    \label{33}
\end{align}
On the other hand, using the assumption that $q\in (d/2, \infty)$ 
and the Sobolev embedding \eqref{emb2} again, we can also get
$$
\nabla u_n\in W^{1,q}(\mR^d)\subseteq W^{\beta,q}(\mR^d)\hookrightarrow L^{2q}(\mR^d) \quad\text{with}\quad \beta=\tfrac{d}{2q}<1.
$$
Thus, for every $\eps>0$, there exists a $C_\eps>0$ such that
\begin{align}
\||\nabla u_n|^2\|_q=\|\nabla u_n\|_{2q}^2\leq \|\nabla u_n\|_{\beta,q}^2\leq \eps\|\nabla^2u_n\|_q^2+C_\eps\|u_n\|_q^2.\label{333}
\end{align}
By the induction assumption, we have
\begin{align*}
&\lambda\|u_{n+1}\|_q+\sqrt{\lambda}\|\nabla u_{n+1}\|_q+\|\nabla^2u_{n+1}\|_{q}\\
&\leq C_d\Big(\|f\|_q+\eps \|G\|_p\|\nabla^2u_n\|_{q}+C_\eps \|G\|_p\|u_n\|_q+\eps\kappa \|\nabla^2u_n\|_{q}^2+\kappa C_\eps \|u_n\|_q^2\Big)\\
&\leq C_d\Big(\|f\|_q+\eps \|G\|_p[C_d\cC_n\|f\|_q]+C_\eps \|G\|_p[\tfrac{C_d\cC_{n}}{\lambda}\|f\|_{q}]\\
&\qquad\qquad+\eps\kappa [C_d\cC_{n}\|f\|_{q}]^2+\kappa C_\eps \big[\tfrac{C_d\cC_{n}}{\lambda}\|f\|_{q}\big]^2\Big)\\
&\leq C_d\|f\|_q\Big(1+\big(\eps C_d\|G\|_p+C_\eps \tfrac{C_d}{\lambda}\|G\|_p+\eps\kappa {C_d}^2\|f\|_{q}+\kappa C_\eps {C_d}^2\tfrac{\|f\|_{q}}{\lambda^2}\big){\cC_n}^2\Big).
\end{align*}
Take
$$
\delta:=\eps C_d\|G\|_p+C_\eps \tfrac{C_d}{\lambda}\|G\|_p+\eps\kappa {C_d}^2\|f\|_{q}+\kappa C_\eps {C_d}^2\tfrac{\|f\|_{q}}{\lambda^2}.
$$
Then by the recursive definition of $\cC_{n+1}$,
we obtain \eqref{indu} holds for every $n\geq 1$. Thus the claim is true.

\vspace{2mm}
Now, for a given $f\in L^q(\mR^d)$, we can first choose $\eps$ small enough so that
$$
\eps C_d\|G\|_p+\eps\kappa {C_d}^2\|f\|_{q}\leq 1/8,
$$
and then take $\lambda_0$ large enough such that
$$
C_\eps \tfrac{C_d}{\lambda_0}\|G\|_p+\kappa C_\eps {C_d}^2\tfrac{\|f\|_{q}}{\lambda_0^2}\leq 1/8
$$
to get by Lemma \ref{nu} that for every $\lambda\geq \lambda_0$,
the sequence $\cC_n$ converges to a finite number $\cC>0$.
This in particular means that
\begin{align}
\sup_{n\geq 1}\Big(\lambda\|u_n\|_q+\sqrt{\lambda}\|\nabla u_n\|_q+\|\nabla^2u_n\|_{q}\Big)\leq C_d\,\cC\|f\|_{q}.\label{un}
\end{align}

Next, we show that $\{u_n\}$ is a Cauchy sequence in $W^{2,q}(\mR^d)$. To this end, note that for $n,m\geq 1$,
$$
u_n(x)-u_m(x)=\int_0^\infty\!\e^{-\lambda t}\cT_t\big(G\cdot[\nabla u_{n-1}-\nabla u_{m-1}]+\kappa[|\nabla u_{n-1}|^2-|\nabla u_{m-1}|^2]\big)(x)\dif t,
$$
we have by H\"older's inequality that
\begin{align*}
\lambda\|u_n-u_m\|_{q}+\|u_n-u_m\|_{2,q}&\leq C_d\Big(\|G\cdot[\nabla u_{n-1}-\nabla u_{m-1}]\|_{q}\\
&\quad+\kappa\big(\|\nabla u_{n-1}\|_{2q}+\|\nabla u_{m-1}\|_{2q}\big)\|\nabla u_{n-1}-\nabla u_{m-1}\|_{2q}\Big).
\end{align*}
Using the same argument as in \eqref{33}--\eqref{333}, we can get that for every $\eps>0$, there exists a $C_\eps>0$ such that
\begin{align*}
\|G\cdot[\nabla u_{n-1}-\nabla u_{m-1}]\|_{q}\leq \eps\|G\|_p\|u_{n-1}-u_{m-1}\|_{2,q}+C_\eps\|G\|_p\|u_n-u_m\|_{q}
\end{align*}
and
\begin{align*}
&\kappa\big(\|\nabla u_{n-1}\|_{2q}+\|\nabla u_{m-1}\|_{2q}\big)\|\nabla u_{n-1}-\nabla u_{m-1}\|_{2q}\\
&\leq \Big(\eps(\|\nabla u_{n-1}\|_{2q}+\|\nabla u_{m-1}\|_{2q})+C_\eps(\|u_n\|_q+\|u_m\|_{q})\Big)\|u_{n-1}-u_{m-1}\|_{2,q}.
\end{align*}
In view of \eqref{un},
we can first take $\eps$ small and then take $\lambda>C_\eps\|G\|_p$ large enough such that
$$
\sup_{n\geq 1}C_d\Big(\eps\|G\|_p+\eps(\|\nabla u_{n-1}\|_{2q}+\|\nabla u_{m-1}\|_{2q})+C_\eps(\|u_n\|_q+\|u_m\|_{q})\Big)<1,
$$
which in turn implies the desired conclusion.

\vspace{2mm}
Based on the above discussions, we know that there exists a limit function $u\in W^{2,q}(\mR^d)$ such that
$$
\lim_{n\rightarrow\infty}\|u_n-u\|_{2,q}=0.
$$
Taking limit on both sides of \eqref{ind}, we get that $u$ satisfies (\ref{inte}), and the estimate \eqref{ess} follows from \eqref{un}. The uniqueness can be proved
using an argument similar to that used in proving the Cauchy property of $u_n$.

\vspace{2mm}
{\bf(ii)}
We only need to notice that in the case  $f\in C_0^\infty(\mR^d)$, the function $u_1$ defined by \eqref{u1} belongs to $C_b^2(\mR^d)$. The proof of the desired results is entirely similar 
to the above by replacing the norm $\|\cdot\|_{2,q}$ with $\|\cdot\|_{C_b^2}$.
We omit the details here.
\end{proof}

\br
From the above proof we can see that if we require the $L^q$-norm of $f$ small enough, then we can get the well-posedness of \eqref{eq2} for every $\lambda>0$.
For large initial data $f$ the solution will blow up. 
Such phenomenon appears naturally in non-linear partial differential equations. Here, we do not require the $L^q$-norm of $f$ small, but we require $\lambda$ large enough instead.
\er

\subsection{Krylov's estimate by the Girsanov transform}
Recall that for every $f\in L^p(\mR^d)$ with $p\in(d\vee2,\infty]$, we have
\begin{align}
\sup_{v\in\mR^d}\mE\exp\left\{\int_0^t|f(v+W_s)|^2\dif s\right\}\leq C\e^{Ct}, 
 \label{aa}
\end{align}
where $C>0$ is a constant depending on $\|f\|_p$.
Using the Girsanov theorem, we can prove
a local Krylov's estimate for the solution of  \eqref{sde1}.
We have the following result.

\bl\label{exis}
Assume that $G\in L^p(\mR^d)$ with $p\in (d\vee2, \infty]$.
Then, for every initial value $y=(x,v)^{\mathrm{T}}\in\mR^{2d}$, there exists a weak solution to \eqref{sde1}. Moreover,
let $(X_t(x),V_t(v))^{\mathrm{T}}$ solve \eqref{sde1} and for every $R>0$, define
\begin{align}
\tau_{R}^v:=\inf\{t\geq 0:|V_t(v)|\geq R\}.    \label{tau}
\end{align}
Then, for every $T>0$ and $f\in L^q(\mR^d)$ with $q\in((d/2)\vee1, \infty]$, we have
\begin{align}
\mE\exp\left\{\int_0^{T\wedge\tau_R^v}\!f(V_s(v))\dif s\right\}\leq C_R\e^{C_RT},    \label{kry}
\end{align}
where $C_R=C(d,x,v,R,\|f\|_{q},\|G\|_p)$ is a positive constant which is uniformly bounded for $(x,v)$ in compact sets.
\el

\begin{proof}
Let $W_t$ be a standard Brownian motion on some probability space $(\Omega,\sF,\mP)$,
and $(\hat X_t,\hat V_t)$ be given by
\begin{equation*}
\left\{ \begin{aligned}
&\dif \hat X_t=\hat V_t\dif t,\quad \hat X_0=x,\\
&\dif \hat V_t=\dif W_t,\quad\,\,\hat V_0=v.
\end{aligned} \right.
\end{equation*}
Define the process
$$
H_t:=W_t+\int_0^t\hat\xi(s)\dif s,
$$
where $\hat\xi(s):=\gamma\hat V_s+\nabla F(\hat X_s)-G(\hat V_s)$,
and set
$$
\hat\Psi_t:=\exp\bigg\{\int_0^t\hat\xi(s)\dif W_s-\frac{1}{2}\int_0^t\hat\xi(s)^2\dif s\bigg\}.
$$
Fix $T>0$ below. For $R>0$, let
$$
\hat\tau_R^v:=\inf\{t\geq 0:|\hat V_t(v)|\geq R\}.
$$
Note that for $t\leq T\wedge\hat\tau_R^v$, we have
$|\hat X_t|\leq x+vT+RT$.
Meanwhile, since $p\in (d\vee2, \infty]$, we have by \eqref{aa} that
\begin{align*}
\sup_{v\in\mR^d}\mE\exp\left\{\int_0^t|G(\hat V_s(v))|^2\dif s\right\}\leq C_1\e^{C_1t}.
\end{align*}
Thus, the process $\{\hat\Psi_{\cdot\wedge\hat\tau^v_R}\}$ is a martingale by Novikov's criterion.
Define a new probability measure by
$$
\frac{\dif \mQ}{\dif \mP}\big|_{\sF_{T\wedge\hat\tau^v_R}}:=\hat\Psi_{T\wedge\hat\tau^v_R}.
$$
Then, by the Girsanov theorem we know $H_t$ is a Brownian motion on the new probability space $(\Omega,\sF,\mQ)$. It also holds that
$$
\dif \hat V_t=-\gamma\hat V_t\dif t-\nabla F(\hat X_t)\dif t+G(\hat V_t)\dif t+\dif H_t,
$$
which means that $(\hat X_t,\hat V_t)$ is a weak solution to \eqref{sde1} until $t\leq\hat\tau_R^v$  under the probability measure $\mQ$.

Now we turn to prove the estimate \eqref{kry}. We have
\begin{align*}
\mE^{\mP}\exp\left\{\int_0^{T\wedge\tau^v_R}\!f(V_s(v))\dif s\right\}&=\mE^{\mQ}\exp\left\{\int_0^{T\wedge\hat\tau^v_R}\!f(\hat V_s(v))\dif s\right\}\\
&=\mE^{\mP}\left(\exp\left\{\int_0^{T\wedge\hat\tau^v_R}\!f(\hat V_s(v))\dif s\right\}\cdot\hat\Psi_{T\wedge\hat\tau^v_R}\right).
\end{align*}
Thus, applying the Cauchy-Schwarz inequality we get
\begin{align*}
\mE^{\mP}\exp\left\{\int_0^{T\wedge\tau^v_R}\!f(V_s(v))\dif s\right\}\leq \left[\mE^{\mP}\exp\bigg\{\int_0^{T\wedge\hat\tau^v_R}\!2f(\hat V_s(v))\dif s\bigg\}\right]^{\frac{1}{2}}\Big(\mE^{\mP}[\hat\Psi_{T\wedge\hat\tau^v_R}^{2}]\Big)^{1/2}.
\end{align*}
Since $q\in ((d/2)\vee1, \infty]$, we can apply \eqref{aa} to $\sqrt{2|f|}$ to get that
$$
\sup_{v\in\mR^d}\mE^{\mP}\exp\bigg\{\int_0^T\!2f(\hat V_s)\dif s\bigg\}\leq C_2\e^{C_2T}.
$$
As for the second term, we can write by the exponential martingale formula that
\begin{align*}
\hat\Psi_{T\wedge\hat\tau^v_R}^{2}=\exp\bigg\{\int_0^{T\wedge\hat\tau^v_R}\!2\xi(s)\dif W_s-\frac{1}{2}\int_0^{T\wedge\hat\tau^v_R}\!|2\xi(s)|^2\dif s+\int_0^{T\wedge\hat\tau^v_R}\!\xi(s)^2\dif s\bigg\},
\end{align*}
whose expectation is less than $C_R\e^{C_RT}$ due to the assumption $G\in L^p(\mR^d)$ 
with $p\in (d\vee2, \infty]$ and \eqref{aa}. 
Combining the above computations, we obtain \eqref{kry} with $C_R$ satisfies the desired properties. The proof is complete.
\end{proof}

\br
By H\"older's inequality, for every $f\in L^q(\mR^d)$ with $q\in ((d/2)\vee1, \infty]$, 
we have
\begin{align*}
&\sup_{v\in\mR^d}\mE\left(\int_0^tf(v+W_s)\dif s\right)=
\sup_{v\in\mR^d}\int^t_0\int_{\mR^d}f(v+y)p(s, y)\dif y\dif s\\
&\le \int^t_0\|p(s, \cdot)\|_p\dif s \|f\|_q \le \int^t_0
(2\pi s)^{-d/(2q)}\dif s \|f\|_q = C\|f\|_q,
\end{align*}
where $\frac1{p}+\frac1{q}=1$ and $C=\frac{2q}{2q-d}(2\pi)^{-\frac{d}{2q}}t^{1-\frac{d}{2q}}$.
By the same argument as above we can see that for every $f\in L^q(\mR^d)$ with $q>(d/2)\vee1$, it holds that
\begin{align}
\mE\left(\int_0^{T\wedge\tau_R^v}\!f(V_s(v))\dif s\right)\leq C_{R}\|f\|_q\e^{C_RT},    \label{kry2}
\end{align}
where $C_{R}$ depends on $\|G\|_p$.
This estimate will be used to prove the pathwise uniqueness of solutions to \eqref{sde1}.
It is obvious that the constant $C_{R}$ will explode as $R\rightarrow\infty$, and the
exponential dependence on $T$ on the right hand side is
mainly caused by the Girsanov transform. Thus, it can not be used to study
the long time behavior of \eqref{sde1}.
\er

\subsection{Krylov's estimate by elliptic equations}
To prove the ergodicity of \eqref{sde1}, we need to strengthen the estimate \eqref{kry2} by showing that the constant on the right hand side does not depend on $R$, and more importantly, depends only linearly on the $t$-variable.
We have the following two results which correspond to the two different assumptions in Theorem \ref{main1}.

\bl\label{lee}
Assume that $G\in L^p(\mR^d)$ with $p\in (d\vee2, \infty)$ 
and that condition {\bf(A)} in Theorem \ref{main1} holds.
Let $(X_t(x),V_t(v))^{\mathrm{T}}$ solve \eqref{sde1}.
Then, for any $f\in L^q(\mR^d)$ with $q\in ((d/2)\vee 1, \infty)$, 
there exists a constant $C=C(d,x,v,\|G\|_p,\|f\|_q)>0$ such that for all $t\ge 0$,
\begin{align}
\mE\left(\int_0^tf(V_s)\dif s\right)\leq C(1+t).    \label{kryt}
\end{align}
\el
\begin{proof}
As in \cite{Kr-Ro,Zh}, we only need to show that \eqref{kryt} holds for every $f\in C_0^\infty(\mR^d)$.
Let $G_n\in C_0^\infty(\mR^d)$ be such that $G_n\rightarrow G$ in $L^p(\mR^d)$, and let $u_n$ be the solution to \eqref{eq2} with $\kappa=1$ and $G$ replaced by $G_n$.
Then by Lemma \ref{ell}, we can apply It\^o's formula to $u_n(V_t)$ and take expectation to get that
\begin{align*}
\mE u_n(V_{t\wedge\tau^v_R})&=u_n(v)+\mE\left(\int_{0}^{t\wedge\tau^v_R}[-\gamma V_s-\nabla F(X_s)]\nabla u_n(V_s)\dif s\right)\\
&\quad+\mE\left(\int_{0}^{t\wedge\tau^v_R}\!\big[\tfrac{1}{2}\Delta u_n+G\cdot\nabla u_n\big](V_s)\dif s\right),
\end{align*}
where $\tau^v_R$ is given by \eqref{tau}.
By the simple inequality $2ab\leq a^2+b^2$, we have
\begin{align*}
\mE u_n(V_{t\wedge\tau^v_R})
&\leq u_n(v)\!+\!\frac{1}{2}\mE\!\left(\int_{0}^{t\wedge\tau^v_R}\![\gamma V_s+\nabla F(X_s)]^2\dif s\right)\!+\!\frac{1}{2}\mE\!\left(\int_{0}^{t\wedge\tau^v_R}\!|G_n-G|^2(V_s)\dif s\right)\\
&\quad+\mE\left(\int_{0}^{t\wedge\tau^v_R}\big[\tfrac{1}{2}\Delta u_n+G_n\cdot\nabla u_n+|\nabla u_n|^2\big](V_s)\dif s\right).
\end{align*}
This in turn yields that
\begin{align*}
\mE\left(\int_{0}^{t\wedge\tau^v_R}\!f(V_s)\dif s\right)&\leq C_\lambda\|u_n\|_\infty(1+t)+\frac{1}{2}\,\mE\left(\int_{0}^{t\wedge\tau^v_R}[\gamma V_s+\nabla F(X_s)]^2\dif s\right)\\
&\quad+\frac{1}{2}\,\mE\left(\int_{0}^{t\wedge\tau^v_R}|G_n-G|^2(V_s)\dif s\right).
\end{align*}
By \eqref{ess}--\eqref{ap}, the assumption $q\in ((d/2)\vee1, \infty)$ 
and the Sobolev embedding \eqref{emb}, we have
$$
\sup_{n\geq 1}\|u_n\|_\infty\leq C_0\|f\|_q.
$$
Thus, thanks to the assumption $G\in L^p(\mR^d)$ with $p\in (d\vee2, \infty)$
and \eqref{kry2}, we can let $n\rightarrow\infty$  to derive that
\begin{align}
\mE\left(\int_0^{t\wedge\tau^v_R}\!f(V_s)\dif s\right)\leq C_1\|f\|_q(1+t)+\frac{1}{2}\mE\left(\int_0^{t\wedge\tau^v_R}\big[|\gamma V_s|+|\nabla F(X_s)|\big]^2\dif s\right),\label{p1}
\end{align}
where $C_1>0$ is a constant.
By a standard approximation argument (see \cite[Lemma 2.4]{H-W}, for instance), 
the above inequality holds for any 
$f\in L^q(\mR^d)$ with $q\in ((d/2)\vee1, \infty)$.

Recall that $\cH(x,v)$ is defined in \eqref{Lya}.
We have by \eqref{lay} and Young's inequality that
\begin{align*}
\mE\cH(X_{t\wedge\tau^v_R},V_{t\wedge\tau^v_R})&\leq \cH(x,v)+c_1t-c_0\mE\left(\int_0^{t\wedge\tau^v_R}\cH(X_s,V_s)\dif s\right)\\
&\qquad+\mE\left(\int_0^{t\wedge\tau^v_R}\<G(V_s),V_s+\tfrac{\gamma}{2}X_s\>\dif s\right)\\
&\leq \cH(x,v)\!+\!c_1t\!-\!c_2\mE\!\left(\int_0^{t\wedge\tau^v_R}\!\cH(X_s,V_s)\dif s\!\right)\!\!+\!c_3\mE\!\left(\int_0^{t\wedge\tau^v_R}\!|G(V_s)|^2\dif s\!\right),
\end{align*}
which in particular implies that
\begin{align}
\mE\left(\int_0^{t\wedge\tau^v_R}\cH(X_s,V_s)\dif s\right)\leq C_2(1+t)+C_3\mE\left(\int_0^{t\wedge\tau^v_R}|G(V_s)|^2\dif s\right), \label{rr}
\end{align}
where $C_2$ depends on $x,v$.
Note that by condition {\bf(A)} and \eqref{lay},
$$
[\gamma v+\nabla F(x)]^2\leq C_4\cH(x,v)
$$
for some constant $C_4>0$. We can take $f=C_3C_4|G|^2$ in \eqref{p1} to get that
\begin{align*}
C_3\mE\left(\int_0^{t\wedge\tau^v_R}|G(V_s)|^2\dif s\right)&\leq C_5\|G\|_p(1+t)+\frac{1}{2C_4}\mE\left(\int_0^{t\wedge\tau^v_R}\big[|\gamma V_s|+|\nabla F(X_s)|\big]^2\dif s\right)\\
&\leq C_5\|G\|_p(1+t)+\frac{1}{2}\mE\left(\int_0^{t\wedge\tau^v_R}\cH(X_s,V_s)\dif s\right).
\end{align*}
This together with \eqref{rr} yields that
\begin{align}
\mE\left(\int_0^{t\wedge\tau^v_R}\cH(X_s,V_s)\dif s\right)\leq C_6\|G\|_p(1+t), \label{aaa}
\end{align}
which in particular implies that the solution does not explode, i.e., $\tau_R^v\rightarrow\infty$ almost everywhere as $R\rightarrow\infty$.
Plugging \eqref{aaa} back into \eqref{p1} and letting $R\rightarrow\infty$, we arrive at
the desired result.
\end{proof}

\bl
Assume that condition {\bf(B)} in Theorem \ref{main1} holds with $p\in (2d, \infty)$
and that $(X_t(x),V_t(v))^{\mathrm{T}}$ solves \eqref{sde1}. Then, 
for every $f\in L^q(\mR^d)$ with $q\in(d, \infty)$,  there exists a 
constant $C=C(d,x,v,\|G\|_p,\|f\|_q)>0$ such that
\begin{align}
\mE\left(\int_0^t\!f(V_s)\dif s\right)\leq C(1+t),\quad\forall t\geq 0.    \label{kryy}
\end{align}
\el

\begin{proof}
We will give the part of the proof that is different
from the corresponding part of the proof of Lemma \ref{lee}.
Let $f\in C_0^\infty(\mR^d)$ and $G_n\in C_0^\infty(\mR^d)$ be such that $G_n\rightarrow G$ in $L^p(\mR^d)$, and let $u_n$ be the solution to \eqref{eq2} with $\kappa=0$ and $G$ replaced by $G_n$.
Then by Lemma \ref{ell}, we can apply It\^o's formula to $u_n(V_t)$ and take expectation to get that
\begin{align*}
\mE u_n(V_{t})&= u_n(v)+\mE\left(\int_0^{t}[-\gamma V_s-\nabla F(X_s)]\nabla u_n(V_s)\dif s\right)\\
&\qquad\qquad+\mE\left(\int_0^{t}\big[\tfrac{1}{2}\Delta u_n+G\cdot\nabla u_n\big](V_s)\dif s\right).
\end{align*}
Now, by \eqref{ess}, \eqref{ap} and the assumption $q\in(d, \infty)$, 
we can use  the Sobolev embedding \eqref{emb} to get that
$$
\sup_{n\geq 1}\big(\|u_n\|_\infty+\|\nabla u_n\|_\infty\big)\leq C_0\|f\|_q.
$$
Thus, we have
\begin{align*}
\mE u_n(V_{t})
&\leq C_\lambda\|u_n\|_\infty(1+t)-\mE\left(\int_0^{t}\!f(V_s)\dif s\right)\\
&\quad+\mE\left(\int_0^{t}\Big(|\gamma V_s|+|\nabla F(X_s)|+|G_n(V_s)-G(V_s)|\Big)\dif s\right)\cdot\|\nabla u_n\|_\infty.
\end{align*}
Letting $n\rightarrow\infty$  we arrive at
\begin{align*}
\mE\left(\int_0^t\!f(V_s)\dif s\right)\leq C_1(1+t)+C_2\mE\left(\int_0^t\big(|\gamma V_s|+|\nabla F(X_s)|\big)\dif s\right).
\end{align*}
Now by condition {\bf(B)} and \eqref{lay} we have
$$
|v|+|\nabla F(x)|\leq C_3\cH(x,v),
$$
which together with (\ref{rr}) yields that
\begin{align*}
\mE\left(\int_0^t\!f(V_s)\dif s\right)\leq C_4(1+t)+C_5\mE\left(\int_0^t|G(V_s)|^2\dif s\right).
\end{align*}
Since $G\in L^p(\mR^{2d})$ with $p\in (2d, \infty)$, 
we can take $f=2C_5|G|^2$ in the above inequality, and the 
desired result follows by the same argument as before.
\end{proof}

\section{Proof of Theorem \ref{main1}}

We will use Zvonkin's transform to remove the singular part of the term $G(V_t)$ 
in \eqref{sde1} when $p\in (d\vee2, \infty)$.
To this end,
let $u$ be the solution to the following equation:
\begin{align}
\lambda u(x)-\tfrac{1}{2}\Delta u(x)-G(x)\cdot\nabla u(x)=G(x).    \label{pde2}
\end{align}
Since $G\in L^p(\mR^d)$ with $p\in (d\vee2, \infty)$, 
we have $u\in W^{2,p}(\mR^d)$. By the Sobolev embedding \eqref{emb}, it holds that $u\in C_b^1(\mR^d)$.
The following result can be proved in a similar way as in \cite{Kr-Ro,Zh} and the only difference
is that we have to use a localization argument to control  the growth of the coefficients.
For convenience, we sketch the proof here.

\bl\label{sde3}
Let $Y_t=(X_t,V_t)^{\mathrm{T}}$ solve \eqref{sde1}
and $\tau_R^v$ be given by \eqref{tau}. Then we also have that for $t\leq \tau_R^v$,
\begin{equation}\label{new}
\left\{ \begin{aligned}
&\dif X_t=V_t\dif t,\\
&\dif V_t=[-\gamma V_t-\nabla F(X_t)]\sigma(V_t)\dif t+[\lambda u(V_t)\dif t-\dif u(V_t)]+\sigma(V_t)\dif W_t,
\end{aligned} \right.
\end{equation}
with initial value $(x,v)^{\mathrm{T}}$, where
$$
\sigma(v):=\big(\mI+\nabla u(v)\big).
$$
\el
\begin{proof}
Let $u_n$ be a mollifying approximation sequence of $u$.
Applying It\^o's formula to $u_n(V_t)$, we get
\begin{align*}
u_n(V_{t\wedge\tau^v_R})&=u_n(v)+\int_0^{t\wedge\tau^v_R}[-\gamma V_s-\nabla F(X_s)]\nabla u_n(V_s)\dif s\\
&\quad+\int_0^{t\wedge\tau^v_R}\big[\tfrac{1}{2}\Delta u_n+G\cdot\nabla u_n\big](V_s)\dif s+\int_0^{t\wedge\tau^v_R}\nabla u_n(V_s)\dif W_s.
\end{align*}
Using Krylov's estimate \eqref{kry2} and letting $n\rightarrow\infty$, we can get that the above equality holds with $u_n$ replaced by $u$. This together with \eqref{pde2} gives a representation that for $t\leq \tau_R^v$,
\begin{align*}
\int_0^tG(V_s)\dif s&=u(v)-u(V_t)+\lambda\int_0^tu(V_s)\dif s\\
&\quad+\int_0^t[-\gamma V_s-\nabla F(X_s)]\nabla u(V_s)\dif s+\int_0^t\nabla u(V_s)\dif W_s.
\end{align*}
Plugging this back into \eqref{sde1} finishes the proof.
\end{proof}

\br
In the new SDE (\ref{new}), the singular part of the drift disappears, 
but the Hamiltonian structure has been totally destroyed. Thus it is difficult to find a Lyapunov function for \eqref{new} (and thus the the original system (\ref{sde1})).
\er

We are now in a position to give:

\begin{proof}[Proof of Theorem \ref{main1}]
We first prove the existence and uniqueness of strong solutions in the case that $G\in L^p(\mR^d)$ with $p\in (d\vee 2, \infty)$.
By Lemma \ref{exis} and the Yamada-Watanabe principle, it suffices to show the pathwise uniqueness for solutions of  \eqref{sde1}. To this end, we consider two solutions $Y_t(y):=\big(X_t(x),V_t(v)\big)^{\mathrm{T}}$ and $Y_t(\hat y):=\big(X_t(\hat x),V_t(\hat v)\big)^{\mathrm{T}}$
of \eqref{sde1}, defined on the same probability space and with respect to the same Brownian motion,
 starting at $y:=(x,v)^{\mathrm{T}}$ and $\hat y:=(\hat x,\hat v)^{\mathrm{T}}$, respectively. Define
$$
\zeta_R:=\inf\{t\geq 0:|V_t(v)|\vee|V_t(\hat v)|\geq R\}.
$$
Let us fix $T>0$ below. We proceed to prove that for every
$\beta\in (0, 1)$, $t\leq T$ and $R>0$,
there exists a constant $C_{R,T}>0$ such that for all $y, \hat y\in \mR^{2d}$ with
$|y|,|\hat y|\leq R$,
\begin{align}\label{esss}
\mE|Y_{t\wedge\zeta_R}(y)-Y_{t\wedge\zeta_R}(\hat y)|^{2\beta}\leq C_{R,T}|y-\hat y|^{2\beta}.
\end{align}
Once this is proved,
we can apply the special case $y=\hat y$ to get the pathwise uniqueness of solutions to
\eqref{sde1}.

Set $\tilde X_t:=X_t(x)-X_t(\hat x)$ and $\tilde V_t:=V_t(v)-V_t(\hat v)$. Then, by Lemma \ref{sde3} we know that the difference process $(\tilde X_t,\tilde V_t)^{\mathrm{T}}$ satisfies the following equation: for $t\leq\zeta_R$,
\begin{equation*}
\left\{ \begin{aligned}
&\dif \tilde X_t=\tilde V_t\dif t,\\
&\dif\tilde V_t=\xi(\tilde X_t,\tilde V_t)\dif t+\Big(\lambda\big[u\big(V_t(v)\big)-u\big(V_t(\hat v)\big)\big]\dif t-\dif [u(V_t(\hat v))-u(V_t(\hat v))]\Big)
\end{aligned} \right.
\end{equation*}
$$
+\big[\sigma\big(V_t(\hat v)\big)-\sigma\big(V_t(\hat v)\big)\big]\dif W_t,
$$
with initial value $\tilde X_0=x-\hat x$ and $\tilde V_0=v-\hat v$, where
\begin{align*}
\xi(\tilde X_t,\tilde V_t)&:=[-\gamma V_t(v)-\nabla F\big(X_t(x)\big)]\sigma\big(V_t(v)\big)-[-\gamma V_t(\hat v)-\nabla F(X_t\big(\hat x)\big)]\sigma\big(V_t(\hat v)\big).
\end{align*}
Note that for $T>0$ and $t\leq T\wedge\zeta_R$,
$$
|X_t(x)|\vee |X_t(\hat x)|\leq R+RT.
$$
Since $u\in W^{2,p}(\mR^d)$, it follows from \eqref{w1p} that there exists a constant $C_{R,T}>0$ such that
\begin{align*}
\xi(\tilde X_{t\wedge\zeta_R},\tilde V_{t\wedge\zeta_R})\leq C_{R,T}(|\tilde V_t|+|\tilde X_t|)+C_{R,T}|\tilde V_t|\cdot\big[g\big(V_t(v)\big)+g\big(V_t(\hat v)\big)\big],
\end{align*}
where the non-negative function $g\in L^p(\mR^d)$ with $p\in (d, \infty)$.
Now, by It\^o's formula, we get for every $t\leq T$,
\begin{align*}
|\tilde X_{t\wedge\zeta_R}|^2+|\tilde V_{t\wedge\zeta_R}|^2&\leq |x-\hat x|^2+C_0|v-\hat v|^2+C_{\lambda,R,T}\int_0^{t\wedge\zeta_R}\!\big(|\tilde X_s|^2+|\tilde V_s|^2\big)\dif s\\
&\quad+C_{R,T}\int_0^{t\wedge\zeta_R}\!|\tilde V_s|^2\dif A_s+\int_0^{t\wedge\zeta_R}\!\big[\sigma\big(V_s(\hat v)\big)-\sigma\big(V_s(\hat v)\big)\big]\tilde V_s\dif W_s,
\end{align*}
where $C_0>0$, and $A_t$ is a continuous increasing process given by
$$
A_t:=\int_0^t\!\Big(g\big(V_s(v)\big)+g\big(V_s(\hat v)\big)\Big)^2\dif s.
$$
In view of \eqref{kry}, we have that for every $\lambda>0$,
$$
\mE\e^{\lambda A_{t\wedge\zeta_R}}\leq 
C_{R,T}\e^{C_{R,T}}<\infty.
$$
This in particular yields (\ref{esss}) by \cite[Lemma 3.8]{XZ2}. 
The case that $G\in L^\infty(\mR^d)$ can be proved by a localization technique as in \cite[Corollary 2.6]{XZ2} or \cite[Theorem 1.3]{Zh1}.

Now we prove the ergodicity of the strong solution. 
In the case when $p\in (d\vee2, \infty)$, 
we can use \eqref{lay}, \eqref{kryt} in the case of {\bf (A)} and \eqref{kryy} in the case of {\bf (B)} to get that
\begin{align*}
\mE\cH(X_t,V_t)&\leq \cH(x,v)+c_1t-c_2\mE\left(\int_0^t\cH(X_s,V_s)\dif s\right)+c_3\mE\left(\int_0^t|G(V_s)|^2\dif s\right)\\
&\leq \cH(x,v)+c_4(1+t)-c_2\mE\left(\int_0^t\cH(X_s,V_s)\dif s\right).
\end{align*}
In the case when $p=\infty$, the inequalities above are obviously valid. 
By Gronwall's inequality, we have
\begin{align}
\sup_{t\geq 0}\mE\cH(X_t,V_t)\leq C<\infty,  \label{unifor}
\end{align}
which implies the existence of invariant distributions.
By Theorem \ref{sf} below, we know that the strong solution is strong Feller and irreducible. Thus, the exponential ergodicity follows by a standard argument (see \cite{Go-Ma}).
\end{proof}

We proceed to show that the unique strong solution is strong Feller and irreducible. Note that \eqref{sde1} is degenerate, the argument used in \cite{Kr-Ro,XZ,XZ2,Zh1} does not apply. Here, we adopt the method developed in \cite{M-S} by making use of  the Girsanov transform.

\bt\label{sf}
The unique strong solution $Y_t=\big(X_t,V_t\big)^{\mathrm{T}}$ is strong Feller as well as irreducible.
\et
\begin{proof}
Consider the following approximation SDE:
\begin{equation*}
\left\{ \begin{aligned}
&\dif X_t^n=V_t^n\dif t,\qquad\qquad\qquad\qquad\qquad\qquad\qquad\quad\qquad\quad\,\, X^n_0=x\in\mR^{d},\\
&\dif V_t^n=-\gamma \chi_n(V^n_t)\dif t-\nabla F_n(X^n_t)\dif t+G(V^n_t)\dif t+\dif W_t,\quad V^n_0=v\in\mR^{d},
\end{aligned} \right.
\end{equation*}
where for $n>0$, $\chi_{n}$ is a smooth  function satisfying
$$
\chi_{n}(x)=\left\{
\begin{aligned}
1 &,\quad  |x|\leq n;\\
0 &,\quad  |x|> 2n,
\end{aligned}
\right.
$$
and $F_n(x):=F(x)\chi_n(x)$.
Then, repeating the  argument as in the proof of \eqref{esss}, we can deduce that for every $(x,v)^{\mathrm{T}}$ and $(\hat x,\hat v)^{\mathrm{T}}\in \mR^{2d}$, there exists a constant $C_{n}>0$ such that
$$
\mE\Big(|X_t^n(x)-X^n_t(\hat x)|+|V^n_t(v)-V^n_t(\hat v)|\Big)\leq C_{n}(|x-\hat x|+|v-\hat v|).
$$
Note that there is no need to localize any more,
because the coefficients $\chi_n$ and $F_n$ are  bounded. This in turn means that
$(X_t^n,V_t^n)^{\mathrm{T}}$ is a
Feller process in the sense that the semigroup of $(X_t^n,V_t^n)^{\mathrm{T}}$ maps bounded continuous functions to bounded continuous functions.
On the other hand, it is obvious that the process
\begin{equation*}
\left\{ \begin{aligned}
&\dif \hat X_t^n=\hat V_t^n\dif t,\quad\,\, X^n_0=x\in\mR^{d},\\
&\dif \hat V_t^n=\dif W_t,\qquad V^n_0=v\in\mR^{d}
\end{aligned} \right.
\end{equation*}
is strong Feller. Thus, by \cite[Theorem 2.1]{M-S}, we have that $Y_t^n:=(X_t^n,V_t^n)^{\mathrm{T}}$ is strong Feller.
That is, for any bounded Borel function $f$ and $t>0$,
\begin{align}
y\mapsto \mE f(Y^n_{t}(y))\mbox{ is continuous.}\label{EB6}
\end{align}
On the other hand, for $R>0$ and $n>R$, define a stopping time
$$
\zeta_{n,R}:=\left\{t\geq 0: \sup_{|y|\leq R}|Y_t(y)|\geq n\right\}.
$$
By the local uniqueness of solutions to \eqref{sde1}, we have that
$$
Y_t(y)=Y^n_t(y),\quad\forall t\in[0,\zeta_{n,R}].
$$
Thus, given $y,\hat y\in \mR^{2d}$ with $|y|,|\hat y|\leq R$, we have
\begin{align*}
|\mE(f(Y_t(y))-f(Y_t(\hat y)))|&\leq\big|\mE\big(f(Y_{t\wedge\zeta_{n,R}}(y))-f(Y_{t\wedge\zeta_{n,R}}(\hat y))1_{t\leq\zeta_{n,R}}\big)\big|\\
&\quad+2\|f\|_\infty \mP(t>\zeta_{n,R}).
\end{align*}
Note that by \eqref{unifor}, $\mP(t>\zeta_{n,R})\rightarrow0$ as $n\rightarrow\infty$, which together with \eqref{EB6} yields the continuity of $y\mapsto \mE(f(Y_t(y)))$.

To show the irreducibility, we have 
by the proof of Lemma  \ref{exis} that
for every $t\in(0,T]$ and open set $A\subseteq\mR^{2d}$,
\begin{align*}
\mP\Big((X_{t\wedge\tau^v_R},V_{t\wedge\tau^v_R})\in A\Big)=\mE\big[\hat \Psi_{t\wedge\hat\tau^v_R}\cdot1_A(\hat Y_{t\wedge\hat\tau^v_R})\big].
\end{align*}
Note that the process $\hat\Psi_t$ is strictly positive. Denote by $\Omega_A:=\{\omega|\hat Y_{t\wedge\hat\tau^v_R}\in A\}$. It holds that for $R$ large enough, $\mP(\Omega_A)>0$ and hence
$$
\mP\Big((X_{t\wedge\tau^v_R},V_{t\wedge\tau^v_R})^{\mathrm{T}}\in A\Big)\geq \int_{\Omega_A}\hat\Psi_{t\wedge\hat\tau^v_R}\dif\mP>0,
$$
which yields the desired result.
\end{proof}

\end{document}